\documentclass[final,3p]{elsarticle}
 \usepackage{graphics}
 \usepackage{graphicx}
 \usepackage{epsfig}
\usepackage{amssymb}
 \usepackage{amsthm}
 \usepackage{lineno}
 \usepackage{amsmath}
   \numberwithin{equation}{section}
\usepackage{mathrsfs}

\newtheorem{thm}{Theorem}[section]

\newtheorem{lem}[thm]{Lemma}

\biboptions{sort&compress,square}
\allowdisplaybreaks
\begin{document}
\begin{frontmatter}
\author{Tong Wu}
\ead{wut977@nenu.edu.cn}
\author{Yong Wang\corref{cor3}}
\ead{wangy581@nenu.edu.cn}
\cortext[cor3]{Corresponding author.}

\address{School of Mathematics and Statistics, Northeast Normal University,
Changchun, 130024, China}
\title{The spectral Einstein functional and the noncommutative residue for manifolds with boundary}
\begin{abstract}
In this paper, we define the spectral Einstein functional associated with the Dirac operator for manifolds with boundary.  And we give the proof of Kastler-Kalau-Walze type theorem for the spectral Einstein functional associated with the Dirac operator on 4-dimensional manifolds with boundary.
\end{abstract}
\end{frontmatter}
\section{Introduction}
 Until now, many geometers have studied noncommutative residues. In \cite{Gu,Wo}, authors found noncommutative residues are of great importance to the study of noncommutative geometry. In \cite{Co1}, Connes used the noncommutative residue to derive a conformal 4-dimensional Polyakov action analogy. Connes showed us that the noncommutative residue on a compact manifold $M$ coincided with the Dixmier's trace on pseudodifferential operators of order $-{\rm {dim}}M$ in \cite{Co2}.
And Connes claimed the noncommutative residue of the square of the inverse of the Dirac operator was proportioned to the Einstein-Hilbert action.  Kastler \cite{Ka} gave a
brute-force proof of this theorem. Kalau and Walze proved this theorem in the normal coordinates system simultaneously in \cite{KW} .
Ackermann proved that
the Wodzicki residue  of the square of the inverse of the Dirac operator ${\rm  Wres}(D^{-2})$ in turn is essentially the second coefficient
of the heat kernel expansion of $D^{2}$ in \cite{Ac}.

On the other hand, Wang generalized the Connes' results to the case of manifolds with boundary in \cite{Wa1,Wa2},
and proved the Kastler-Kalau-Walze type theorem for the Dirac operator and the signature operator on lower-dimensional manifolds
with boundary \cite{Wa3}. In \cite{Wa3,Wa4}, Wang computed $\widetilde{{\rm Wres}}[\pi^+D^{-1}\circ\pi^+D^{-1}]$ and $\widetilde{{\rm Wres}}[\pi^+D^{-2}\circ\pi^+D^{-2}]$, where the two operators are symmetric, in these cases the boundary term vanished. But for $\widetilde{{\rm Wres}}[\pi^+D^{-1}\circ\pi^+D^{-3}]$, Wang got a nonvanishing boundary term \cite{Wa5}, and give a theoretical explanation for gravitational action on boundary. In others words, Wang provides a kind of method to study the Kastler-Kalau-Walze type theorem for manifolds with boundary. In \cite{DL}, the authors defined bilinear functionals of vector fields and differential forms, the densities of which yield the  metric and Einstein tensors on even-dimensional Riemannian manifolds. Motivated by \cite{DL}, we define  the spectral Einstein functional associated with the Dirac operator for manifolds with boundary, and the motivation of this paper is
to compute the noncommutative residue $\widetilde{{\rm Wres}}[\pi^+(\nabla_X^{S(TM)}\nabla_Y^{S(TM)}D^{-2})\circ\pi^+(D^{-2})]$ on 4-dimensional compact manifolds. \\
\indent The paper is organized in the following way. In Section \ref{section:2}, we define the spectral Einstein functional associated with the Dirac operator and and get the noncommutative residue for manifolds without boundary. In Section \ref{section:3}, we prove the Kastler-Kalau-Walze type theorem for the spectral Einstein functional associated with the Dirac operator on 4-dimensional manifolds with boundary.
\section{The spectral Einstein functional associated with the Dirac operator}
\label{section:2}
Firstly we recall the definition of Dirac operator. Let $M$ be an $n$-dimensional ($n\geq 3$) oriented compact Riemannian manifold with a Riemannian metric $g^{M}$ and let $\nabla^L$ be the Levi-Civita connection about $g^{M}$. In the fixed orthonormal frame $\{e_1,\cdots,e_n\}$, the connection matrix $(\omega_{s,t})$ is defined by
\begin{equation}
\label{a2}
\nabla^L(e_1,\cdots,e_n)= (e_1,\cdots,e_n)(\omega_{s,t}).
\end{equation}
\indent Let $\epsilon (e_j^*)$,~$\iota (e_j^*)$ be the exterior and interior multiplications respectively, where $e_j^*=g^{TM}(e_j,\cdot)$.
Write
\begin{equation}
\label{a3}
\widehat{c}(e_j)=\epsilon (e_j^* )+\iota
(e_j^*);~~
c(e_j)=\epsilon (e_j^* )-\iota (e_j^* ),
\end{equation}
which satisfies
\begin{align}
\label{a4}
&\widehat{c}(e_i)\widehat{c}(e_j)+\widehat{c}(e_j)\widehat{c}(e_i)=2g^{M}(e_i,e_j);~~\nonumber\\
&c(e_i)c(e_j)+c(e_j)c(e_i)=-2g^{M}(e_i,e_j);~~\nonumber\\
&c(e_i)\widehat{c}(e_j)+\widehat{c}(e_j)c(e_i)=0.\nonumber\\
\end{align}
By \cite{Y}, we have the Dirac operator
\begin{align}
\label{a5}
&D=\sum^n_{i=1}c(e_i)[e_i-\frac{1}{4}\sum_{s,t}\omega_{s,t}
(e_i)c(e_s)c(e_t)].\nonumber\\
\end{align}
\indent We define $\nabla_X^{S(TM)}:=X+\frac{1}{4}\Sigma_{ij}\langle\nabla_X^L{e_i},e_j\rangle c(e_i)c(e_j)$ is a spin connection. Set $A(X)=\frac{1}{4}\Sigma_{ij}\langle\nabla_X^L{e_i},e_j\rangle c(e_i)c(e_j)$, then
\begin{align}\label{ddd}
\nabla_X^{S(TM)}\nabla_Y^{S(TM)}&=[X+A(X)][Y+A(Y)]\nonumber\\
&=XY+X\cdot A(Y)+A(X)Y+A(X)A(Y)\nonumber\\
&=XY+X[A(Y)]+A(Y)X+A(X)Y+A(X)A(Y),\nonumber\\
\end{align}
where $X=\Sigma_{j=1}^nX_j\partial_{x_j}, Y=\Sigma_{l=1}^nY_l\partial_{x_l}$ and $XY=\Sigma_{j,l=1}^n(X_jY_l\partial_{x_j}\partial_{x_l}+X_j\frac{\partial_{Y_l}}{\partial_{x_j}}\partial_{x_l}).$ \\
\indent Let $g^{ij}=g(dx_{i},dx_{j})$, $\xi=\sum_{k}\xi_{j}dx_{j}$ and $\nabla^L_{\partial_{i}}\partial_{j}=\sum_{k}\Gamma_{ij}^{k}\partial_{k}$,  we denote that
\begin{align}
&\sigma_{i}=-\frac{1}{4}\sum_{s,t}\omega_{s,t}
(e_i)c(e_s)c(e_t)
;~~~\xi^{j}=g^{ij}\xi_{i};~~~~\Gamma^{k}=g^{ij}\Gamma_{ij}^{k};~~~~\sigma^{j}=g^{ij}\sigma_{i};
~~~~a^{j}=g^{ij}a_{i}.
\end{align}
And by $\partial_{x_j}=-\sqrt{-1}\xi_j$, we have the following lemmas.
\begin{lem}\label{lem3} The following identities hold:
\begin{align}
\label{b22}
&\sigma_{0}(\nabla_X^{S(TM)}\nabla_Y^{S(TM)})=X[A(Y)]+A(X)A(Y);\nonumber\\
&\sigma_{1}(\nabla_X^{S(TM)}\nabla_Y^{S(TM)})=\sqrt{-1}\Sigma_{j,l=1}^nX_j\frac{\partial_{Y_l}}{\partial_{x_j}}\partial_{x_l}+\sqrt{-1}\Sigma_jA(Y)X_j\xi_j+\sqrt{-1}\Sigma_lA(Y)Y_l\xi_l;\nonumber\\
&\sigma_{2}(\nabla_X^{S(TM)}\nabla_Y^{S(TM)})=-\Sigma_{j,l=1}^nX_jY_l\xi_j\xi_l.\nonumber\\
\end{align}
\end{lem}
\begin{lem}\cite{Ka}\label{lem356} The following identities hold:
\begin{align}
\label{b22222}
&\sigma_{-2}(D^{-2})=|\xi|^{-2};\nonumber\\
&\sigma_{-3}(D^{-2})=-\sqrt{-1}|\xi|^{-4}\xi_k(\Gamma^k-2\delta^k)-\sqrt{-1}|\xi|^{-6}2\xi^j\xi_\alpha\xi_\beta\partial_jg^{\alpha\beta}.\nonumber\\
\end{align}
\end{lem}
\indent Write
 \begin{eqnarray}
D_x^{\alpha}&=(-i)^{|\alpha|}\partial_x^{\alpha};
~\sigma(D_t)=p_1+p_0;
~(\sigma(D_t)^{-1})=\sum^{\infty}_{j=1}q_{-j}.
\end{eqnarray}

\indent By the composition formula of pseudodifferential operators, we have
\begin{align}
1=\sigma(D\circ D^{-1})&=\sum_{\alpha}\frac{1}{\alpha!}\partial^{\alpha}_{\xi}[\sigma(D)]
D_x^{\alpha}[\sigma(D^{-1})]\nonumber\\
&=(p_1+p_0)(q_{-1}+q_{-2}+q_{-3}+\cdots)\nonumber\\
&~~~+\sum_j(\partial_{\xi_j}p_1+\partial_{\xi_j}p_0)(
D_{x_j}q_{-1}+D_{x_j}q_{-2}+D_{x_j}q_{-3}+\cdots)\nonumber\\
&=p_1q_{-1}+(p_1q_{-2}+p_0q_{-1}+\sum_j\partial_{\xi_j}p_1D_{x_j}q_{-1})+\cdots,
\end{align}
so
\begin{equation}
q_{-1}=p_1^{-1};~q_{-2}=-p_1^{-1}[p_0p_1^{-1}+\sum_j\partial_{\xi_j}p_1D_{x_j}(p_1^{-1})].
\end{equation}
Then,
\begin{align}\label{mki}
\sigma_{0}(\nabla_X^{S(TM)}\nabla_Y^{S(TM)}D^{-2})&=-\Sigma_{j,l=1}^nX_jY_l\xi_j\xi_l|\xi|^{-2}.
\end{align}
By \cite{DL}, we have the following theorem
\begin{thm}\label{thm2} If $M$ is an $n$-dimensional compact oriented manifolds without boundary, and $n$ is even, then we get the following equality:
\begin{align}
\label{a29}
{\rm Wres}[\nabla_X^{S(TM)}\nabla_Y^{S(TM)}D^{-n}]
&=\frac{(2\pi)^{\frac{n}{2}}}{3(\frac{n}{2}-1)!}\int_{M}[Ric(X,Y)-\frac{1}{2}sg(X,Y)]d{\rm Vol_{M}}+\frac{(2\pi)^{\frac{n}{2}}}{4(\frac{n}{2}-1)!}\int_{M}sg(X,Y)d{\rm Vol_{M}},
\end{align}
where $s$ is the scalar curvature, ${\rm Vol_{M}}$ is the volume of $M$ and $Ric$ denotes Ricci tensor on $M$.
\end{thm}

\section{The noncommutative residue for $4$-dimensional manifolds with boundary}
\label{section:3}
 In this section, we compute the noncommutative residue $\widetilde{{\rm Wres}}[\pi^+(\nabla_X^{S(TM)}\nabla_Y^{S(TM)}D^{-2})\circ\pi^+(D^{-2})]$ on $4$-dimensional oriented compact manifolds with boundary. We firstly recall some basic facts and formulas about Boutet de
Monvel's calculus and the definition of the noncommutative residue for manifolds with boundary which will be used in the following. For more details, see Section 2 in \cite{Wa3}.\\
 \indent Let $M$ be a 4-dimensional compact oriented manifold with boundary $\partial M$.
We assume that the metric $g^{M}$ on $M$ has the following form near the boundary,
\begin{equation}
\label{b1}
g^{M}=\frac{1}{h(x_{n})}g^{\partial M}+dx _{n}^{2},
\end{equation}
where $g^{\partial M}$ is the metric on $\partial M$ and $h(x_n)\in C^{\infty}([0, 1)):=\{\widehat{h}|_{[0,1)}|\widehat{h}\in C^{\infty}((-\varepsilon,1))\}$ for
some $\varepsilon>0$ and $h(x_n)$ satisfies $h(x_n)>0$, $h(0)=1,$ where $x_n$ denotes the normal directional coordinate.\\
\indent  Then similar to \cite{Wa3}, we can compute the noncommutative residue
\begin{align}
\label{b14}
&\widetilde{{\rm Wres}}[\pi^+(\nabla_X^{S(TM)}\nabla_Y^{S(TM)}D^{-2})\circ\pi^+(D^{-2})]\nonumber\\
&=\int_M\int_{|\xi|=1}{\rm
trace}_{\wedge^*T^*M\bigotimes\mathbb{C}}[\sigma_{-4}(\nabla_X^{S(TM)}\nabla_Y^{S(TM)}D^{-2}\circ D^{-2})]\sigma(\xi)dx+\int_{\partial M}\Phi,
\end{align}
where
\begin{align}
\label{b15}
\Phi &=\int_{|\xi'|=1}\int^{+\infty}_{-\infty}\sum^{\infty}_{j, k=0}\sum\frac{(-i)^{|\alpha|+j+k+1}}{\alpha!(j+k+1)!}
\times {\rm trace}_{\wedge^*T^*M\bigotimes\mathbb{C}}[\partial^j_{x_n}\partial^\alpha_{\xi'}\partial^k_{\xi_n}\sigma^+_{r}(\nabla_X^{S(TM)}\nabla_Y^{S(TM)}D^{-2})(x',0,\xi',\xi_n)
\nonumber\\
&\times\partial^\alpha_{x'}\partial^{j+1}_{\xi_n}\partial^k_{x_n}\sigma_{l}(D^{-2})(x',0,\xi',\xi_n)]d\xi_n\sigma(\xi')dx',
\end{align}
and the sum is taken over $r+l-k-j-|\alpha|=-3,~~r\leq 0,~~l\leq-2$.\\

\indent By Theorem \ref{thm2}, we can compute the interior of $\widetilde{{\rm Wres}}[\pi^+(\nabla_X^{S(TM)}\nabla_Y^{S(TM)}D^{-2})\circ\pi^+(D^{-2})]$,\\
 we get
\begin{align}
\label{a16}
{\rm Wres}[\nabla_X^{S(TM)}\nabla_Y^{S(TM)}D^{-4}]
&=\frac{4\pi^2}{3}\int_{M}[Ric(X,Y)-\frac{1}{2}sg(X,Y)]d{\rm Vol_{M}}+\pi^2\int_{M}sg(X,Y)d{\rm Vol_{M}},
\end{align}
\indent Now we  need to compute $\int_{\partial M} \Phi$. When $n=4$, then ${\rm tr}_{S(TM)}[{\rm \texttt{id}}]={\rm dim}(\wedge^*(\mathbb{R}^2))=4$, the sum is taken over $
r+l-k-j-|\alpha|=-3,~~r\leq 0,~~l\leq-2,$ then we have the following five cases:
~\\
\noindent  {\bf case a)~I)}~$r=0,~l=-2,~k=j=0,~|\alpha|=1$.\\
\noindent By (\ref{b15}), we get
\begin{equation}
\label{b24}
\Phi_1=-\int_{|\xi'|=1}\int^{+\infty}_{-\infty}\sum_{|\alpha|=1}
 {\rm tr}[\partial^\alpha_{\xi'}\pi^+_{\xi_n}\sigma_{0}(\nabla_X^{S(TM)}\nabla_Y^{S(TM)}D^{-2})\times
 \partial^\alpha_{x'}\partial_{\xi_n}\sigma_{-2}(D^{-2})](x_0)d\xi_n\sigma(\xi')dx'.
\end{equation}
By Lemma 2.2 in \cite{Wa3}, for $i<n$, then
\begin{equation}
\label{b25}
\partial_{x_i}\sigma_{-2}({D}^{-2})(x_0)=
\partial_{x_i}(|\xi|^{-2})(x_0)=
-\frac{\partial_{x_i}(|\xi|^{2})(x_0)}{|\xi|^4}=0,
\end{equation}
\noindent so $\Phi_1=0$.\\
 \noindent  {\bf case a)~II)}~$r=0,~l=-2,~k=|\alpha|=0,~j=1$.\\
\noindent By (\ref{b15}), we get
\begin{equation}
\label{b26}
\Phi_2=-\frac{1}{2}\int_{|\xi'|=1}\int^{+\infty}_{-\infty} {\rm
trace} [\partial_{x_n}\pi^+_{\xi_n}\sigma_{0}(\nabla_X^{S(TM)}\nabla_Y^{S(TM)}D^{-2})\times
\partial_{\xi_n}^2\sigma_{-2}(D^{-2})](x_0)d\xi_n\sigma(\xi')dx'.
\end{equation}
\noindent By Lemma \ref{lem3}, we have\\
\begin{eqnarray}\label{b237}
\partial_{\xi_n}^2\sigma_{-2}((D^{-2}))(x_0)=\partial_{\xi_n}^2(|\xi|^{-2})(x_0)=\frac{6\xi_n^2-2}{(1+\xi_n^2)^3}.
\end{eqnarray}
\begin{eqnarray}\label{b27}
\partial_{x_n}\sigma_{0}(\nabla_X^{S(TM)}\nabla_Y^{S(TM)}D^{0})=\partial_{x_n}(-\Sigma_{j,l=1}^nX_jY_l\xi_j\xi_l|\xi|^{-2})=\frac{\Sigma_{j,l=1}^nX_jY_l\xi_j\xi_lh'(0)|\xi'|^2}{(1+\xi_n^2)^2}.
\end{eqnarray}
Then, we have
\begin{align}\label{b28}
\pi^+_{\xi_n}\partial_{x_n}\sigma_{0}(\nabla_X^{S(TM)}\nabla_Y^{S(TM)}D^{-2})&=\partial_{x_n}\pi^+_{\xi_n}\sigma_{0}(\nabla_X^{S(TM)}\nabla_Y^{S(TM)}D^{-2})\nonumber\\
&=-\frac{i\xi_n}{4(\xi_n-i)^2}\Sigma_{j,l=1}^{n-1}X_jY_l\xi_j\xi_lh'(0)+\frac{2-i\xi_n}{4(\xi_n-i)^2}X_nY_nh'(0)\nonumber\\
&-\frac{i}{4(\xi_n-i)^2}\Sigma_{j=1}^{n-1}X_jY_n\xi_j-\frac{i}{4(\xi_n-i)^2}\Sigma_{l=1}^{n-1}X_nY_l\xi_l.\nonumber\\
\end{align}
Moreover,
\begin{align}\label{33}
&{\rm
tr} [\partial_{x_n}\pi^+_{\xi_n}\sigma_{0}(\nabla_X^{S(TM)}\nabla_Y^{S(TM)}D^{-2})\times
\partial_{\xi_n}^2\sigma_{-2}(D^{-2})](x_0)\nonumber\\
&=2\frac{1+\xi_ni-3\xi_n^3i-i}{(\xi_n-i)^5(\xi_n+i)^3}\Sigma_{j,l=1}^{n-1}X_jY_l\xi_j\xi_lh'(0)+2\frac{1+\xi_ni-3\xi_n^3i-i}{(\xi_n-i)^5(\xi_n+i)^3}X_nY_nh'(0)\nonumber\\
&+2\frac{(1-3\xi_n^2)i}{(\xi_n-i)^5(\xi_n+i)^3}\Sigma_{j=1}^{n-1}X_jY_n\xi_j+2\frac{(1-3\xi_n^2)i}{(\xi_n-i)^5(\xi_n+i)^3}\Sigma_{l=1}^{n-1}X_nY_l\xi_l.\nonumber\\
\end{align}
Therefore, we get
\begin{align}\label{35}
\Phi_2&=-\frac{1}{2}\int_{|\xi'|=1}\int^{+\infty}_{-\infty}\bigg\{2\frac{1+\xi_ni-3\xi_n^3i-i}{(\xi_n-i)^5(\xi_n+i)^3}\Sigma_{j,l=1}^{n-1}X_jY_l\xi_j\xi_lh'(0)+2\frac{1+\xi_ni-3\xi_n^3i-i}{(\xi_n-i)^5(\xi_n+i)^3}X_nY_nh'(0)\bigg\}d\xi_n\sigma(\xi')dx'\nonumber\\
&=-\Sigma_{j,l=1}^{n-1}X_jY_lh'(0)\Omega_3\int_{\Gamma^{+}}\frac{1+\xi_ni-3\xi_n^3i-i}{(\xi_n-i)^5(\xi_n+i)^3}\xi_j\xi_ld\xi_{n}dx'-X_nY_nh'(0)\Omega_3\int_{\Gamma^{+}}\frac{1+\xi_ni-3\xi_n^3i-i}{(\xi_n-i)^5(\xi_n+i)^3}d\xi_{n}dx'\nonumber\\
&=-\Sigma_{j,l=1}^{n-1}X_jY_lh'(0)\Omega_3\frac{2\pi i}{4!}\left[\frac{1+\xi_ni-3\xi_n^3i-i}{(\xi_n+i)^3}\right]^{(4)}\bigg|_{\xi_n=i}dx'\nonumber\\
&-X_nY_nh'(0)\Omega_3\frac{2\pi i}{4!}\left[\frac{1+\xi_ni-3\xi_n^3i-i}{(\xi_n+i)^3}\right]^{(4)}\bigg|_{\xi_n=i}dx'\nonumber\\
&=\left(\frac{13\pi}{24}\Sigma_{j=1}^{n-1}X_jY_j+\frac{13}{32}X_nY_n\right)h'(0)\pi\Omega_3dx',
\end{align}
where ${\rm \Omega_{3}}$ is the canonical volume of $S^{3}.$\\
\noindent  {\bf case a)~III)}~$r=0,~l=-2,~j=|\alpha|=0,~k=1$.\\
\noindent By (\ref{b15}), we get
\begin{align}\label{36}
\Phi_3&=-\frac{1}{2}\int_{|\xi'|=1}\int^{+\infty}_{-\infty}
{\rm trace} [\partial_{\xi_n}\pi^+_{\xi_n}\sigma_{0}(\nabla_X^{S(TM)}\nabla_Y^{S(TM)}D^{-2})\times
\partial_{\xi_n}\partial_{x_n}\sigma_{-2}(D^{-2})](x_0)d\xi_n\sigma(\xi')dx'\nonumber\\
&=\frac{1}{2}\int_{|\xi'|=1}\int^{+\infty}_{-\infty}
{\rm trace} [\partial_{\xi_n}^2\pi^+_{\xi_n}\sigma_{0}(\nabla_X^{S(TM)}\nabla_Y^{S(TM)}D^{-2})\times
\partial_{x_n}\sigma_{-2}(D^{-2})](x_0)d\xi_n\sigma(\xi')dx'.
\end{align}
\noindent By \cite{Ka}, we have
\begin{eqnarray}\label{37}
\partial_{x_n}\sigma_{-2}(D^{-2})(x_0)|_{|\xi'|=1}
=-\frac{h'(0)}{(1+\xi_n^2)^2}.
\end{eqnarray}
\begin{align}\label{38}
\pi^+_{\xi_n}\sigma_{0}(\nabla_X^{S(TM)}\nabla_Y^{S(TM)}D^{-2})&=\frac{i}{2(\xi_n-i)}\Sigma_{j,l=1}^{n-1}X_jY_l\xi_j\xi_l-\frac{1}{2(\xi_n-i)}X_nY_n-\frac{1}{2(\xi_n-i)}\Sigma_{j=1}^{n-1}X_jY_n\xi_j\nonumber\\
&-\frac{1}{2(\xi_n-i)}\Sigma_{l=1}^{n-1}X_nY_l\xi_l.
\end{align}
Then, we have
\begin{align}\label{mmmmm}
\partial_{\xi_n}^2\pi^+_{\xi_n}\sigma_{0}(\nabla_X^{S(TM)}\nabla_Y^{S(TM)}D^{-2})=\frac{i}{(\xi_n-i)^3}\Sigma_{j,l=1}^{n-1}X_jY_l\xi_j\xi_l-\frac{1}{(\xi_n-i)^3}X_nY_n.
\end{align}
\begin{align}\label{39}
&{\rm tr} [\partial_{\xi_n}\pi^+_{\xi_n}\sigma_{0}(\nabla_X^{S(TM)}\nabla_Y^{S(TM)}D^{-2})\times
\partial_{\xi_n}\partial_{x_n}\sigma_{-2}(D^{-2})](x_0)\nonumber\\
&=-4\frac{h'(0)i}{(\xi_n-i)^5(\xi_n+i)^2}\Sigma_{j,l=1}^{n-1}X_jY_l\xi_j\xi_l+4\frac{h'(0)}{(\xi_n-i)^5(\xi_n+i)^2}X_nY_n.\nonumber\\
\end{align}
Therefore, we get
\begin{align}\label{41}
\Phi_3&=\frac{1}{2}\int_{|\xi'|=1}\int^{+\infty}_{-\infty}
\bigg(-4\frac{h'(0)i}{(\xi_n-i)^5(\xi_n+i)^2}\Sigma_{j,l=1}^{n-1}X_jY_l\xi_j\xi_l+4\frac{h'(0)}{(\xi_n-i)^5(\xi_n+i)^2}X_nY_n\bigg)d\xi_n\sigma(\xi')dx'\nonumber\\
&=-2\Sigma_{j,l=1}^{n-1}X_jY_lh'(0)\Omega_3\int_{\Gamma^{+}}\frac{i}{(\xi_n-i)^5(\xi_n+i)^2}\xi_j\xi_ld\xi_{n}dx'+2X_nY_nh'(0)\Omega_3\int_{\Gamma^{+}}\frac{1}{(\xi_n-i)^5(\xi_n+i)^2}d\xi_{n}dx'\nonumber\\
&=-2\Sigma_{j,l=1}^{n-1}X_jY_lh'(0)\Omega_3\frac{2\pi i}{4!}\left[\frac{i}{(\xi_n+i)^2}\right]^{(4)}\bigg|_{\xi_n=i}dx'+2X_nY_nh'(0)\Omega_3\frac{2\pi i}{4!}\left[\frac{1}{(\xi_n+i)^2}\right]^{(4)}\bigg|_{\xi_n=i}dx'\nonumber\\
&=\left(\frac{5\pi}{12}\Sigma_{j=1}^{n-1}X_jY_j+\frac{5i}{16}X_nY_n\right)h'(0)\pi\Omega_3dx'.
\end{align}
\noindent  {\bf case b)}~$r=0,~l=-3,~k=j=|\alpha|=0$.\\
\noindent By (\ref{b15}), we get
\begin{align}\label{42}
\Phi_4&=-i\int_{|\xi'|=1}\int^{+\infty}_{-\infty}{\rm trace} [\pi^+_{\xi_n}\sigma_{0}(\nabla_X^{S(TM)}\nabla_Y^{S(TM)}D^{-2})\times
\partial_{\xi_n}\sigma_{-3}(D^{-2})](x_0)d\xi_n\sigma(\xi')dx'\nonumber\\
&=i\int_{|\xi'|=1}\int^{+\infty}_{-\infty}{\rm trace} [\partial_{\xi_n}\pi^+_{\xi_n}\sigma_{0}(\nabla_X^{S(TM)}\nabla_Y^{S(TM)}D^{-2})\times
\sigma_{-3}(D^{-2})](x_0)d\xi_n\sigma(\xi')dx'.
\end{align}
 By Lemma \ref{lem3}, we have
\begin{align}\label{43}
\sigma_{-3}(D^{-2})(x_0)|_{|\xi'|=1}=-\frac{i}{(1+\xi_n^2)^2}\left(-\frac{1}{2}h'(0)\sum_{k<n}\xi_nc(e_k)c(e_n)+\frac{5}{2}h'(0)\xi_n\right)-\frac{2ih'(0)\xi_n}{(1+\xi_n^2)^3}.
\end{align}
\begin{align}\label{45}
\partial_{\xi_n}\pi^+_{\xi_n}\sigma_{0}(\nabla_X^{S(TM)}\nabla_Y^{S(TM)}D^{-2})&=-\frac{i}{2(\xi_n-i)^2}\Sigma_{j,l=1}^{n-1}X_jY_l\xi_j\xi_l-\frac{1}{2(\xi_n-i)^2}X_nY_n\nonumber\\
&+\frac{1}{2(\xi_n-i)^2}\Sigma_{j=1}^{n-1}X_jY_n\xi_j+\frac{1}{2(\xi_n-i)^2}\Sigma_{l=1}^{n-1}X_nY_l\xi_l.
\end{align}
Then, we have
\begin{align}\label{39}
&{\rm tr}[\partial_{\xi_n}\pi^+_{\xi_n}\sigma_{0}(\nabla_X^{S(TM)}\nabla_Y^{S(TM)}D^{-2})\times
\sigma_{-3}(D^{-2})](x_0)\nonumber\\
&=-\frac{h'(0)(5\xi_n^2-5+4\xi_n)} {(\xi_n-i)^5(\xi_n+i)^3}\Sigma_{j,l=1}^{n-1}X_jY_l\xi_j\xi_l+\frac{h'(0)i(5\xi_n^3-\xi_n)}{(\xi_n-i)^5(\xi_n+i)^3}X_nY_n.\nonumber\\
\end{align}
Therefore, we get
\begin{align}\label{41}
\Phi_4&=i\int_{|\xi'|=1}\int^{+\infty}_{-\infty}
\bigg(-\frac{h'(0)(5\xi_n^2-5+4\xi_n)} {(\xi_n-i)^5(\xi_n+i)^3}\Sigma_{j,l=1}^{n-1}X_jY_l\xi_j\xi_l+\frac{h'(0)i(5\xi_n^3-\xi_n)}{(\xi_n-i)^5(\xi_n+i)^3}X_nY_n\bigg)d\xi_n\sigma(\xi')dx'\nonumber\\
&=-i\Sigma_{j,l=1}^{n-1}X_jY_lh'(0)\Omega_3\int_{\Gamma^{+}}\frac{5\xi_n^2-5+4\xi_n}{(\xi_n-i)^5(\xi_n+i)^2}\xi_j\xi_ld\xi_{n}dx'+iX_nY_nh'(0)\Omega_3\int_{\Gamma^{+}}\frac{5\xi_n^3-\xi_n}{(\xi_n-i)^5(\xi_n+i)^2}d\xi_{n}dx'\nonumber\\
&=-i\Sigma_{j,l=1}^{n-1}X_jY_lh'(0)\Omega_3\frac{2\pi i}{4!}\left[\frac{5\xi_n^2-5+4\xi_n}{(\xi_n+i)^2}\right]^{(4)}\bigg|_{\xi_n=i}dx'+iX_nY_nh'(0)\Omega_3\frac{2\pi i}{4!}\left[\frac{5\xi_n^3-\xi_n}{(\xi_n+i)^2}\right]^{(4)}\bigg|_{\xi_n=i}dx'\nonumber\\
&=\left(\frac{(1-5i)\pi}{12}\Sigma_{j=1}^{n-1}X_jY_j+\frac{11i}{16}X_nY_n\right)h'(0)\pi\Omega_3dx'.
\end{align}
\noindent {\bf  case c)}~$r=-1,~\ell=-2,~k=j=|\alpha|=0$.\\
By (\ref{b15}), we get
\begin{align}\label{61}
\Phi_5=-i\int_{|\xi'|=1}\int^{+\infty}_{-\infty}{\rm trace} [\pi^+_{\xi_n}\sigma_{-1}(\nabla_X^{S(TM)}\nabla_Y^{S(TM)}D^{-2})\times
\partial_{\xi_n}\sigma_{-2}(D^{-2})](x_0)d\xi_n\sigma(\xi')dx'.
\end{align}
By Lemma \ref{lem3}, we have
\begin{align}\label{62}
\partial_{\xi_n}\sigma_{-2}(D^{-2})|_{|\xi'|=1}=-\frac{2\xi_n}{(\xi_n^2+1)^2}.
\end{align}
Since,
\begin{align}\label{63}
\sigma_{-1}(\nabla_X^{S(TM)}\nabla_Y^{S(TM)}D^{-2})&=i\sum_jA(Y)X_j\xi_j|\xi|^{-2}+i\sum_lA(X)Y_l\xi_l|\xi|^{-2}+i\sum_{j,l=1}^{n-1}X_j\frac{\partial_{Y_l}}{\partial_{X_j}}\xi_l|\xi|^{-2}\nonumber\\
&+i\sum_{j,l=1}^{n-1}X_jY_l\xi_j\xi_l|\xi|^{-4}\xi_k(\Gamma^k-2\delta^k)+2i\sum_{j,l=1}^{n-1}X_jY_l\xi_j\xi_l|\xi|^{-6}\xi^j\xi_\alpha\xi_\beta\partial_jg^{\alpha\beta}.\nonumber\\
\end{align}
By (\ref{62}) and (\ref{63}), we have
\begin{align}\label{71}
&{\rm tr}[\pi^+_{\xi_n}\sigma_{-1}(\nabla_X^{S(TM)}\nabla_Y^{S(TM)}D^{-2})\times
\partial_{\xi_n}\sigma_{-2}(D^{-2})](x_0)|_{|\xi'|=1}\nonumber\\
&=\frac{h'(0)(2\xi_n^2-\xi_n-2\xi_ni)}{(\xi_n-i)^4(\xi_n+i)^2}\sum_{j,l=1}^{n-1}X_jY_l\xi_j\xi_l\sum_{k<n}\xi_kc(e_k)c(e_n)+\frac{h'(0)(17\xi_ni-\xi_n^2+4\xi_n^3i)}{(\xi_n-i)^5(\xi_n+i)^2}\sum_{j,l=1}^{n-1}X_jY_l\xi_j\xi_l.\nonumber\\
\end{align}
So we have
\begin{align}\label{74}
\Phi_5&=-i\int_{|\xi'|=1}\int^{+\infty}_{-\infty}{\rm tr}[\pi^+_{\xi_n}\sigma_{-1}(LD^{-1})\times
\partial_{\xi_n}\sigma_{-2}(D^{-1})](x_0)d\xi_n\sigma(\xi')dx'\nonumber\\
&=-i\int_{|\xi'|=1}\int^{+\infty}_{-\infty}\frac{h'(0)(17\xi_ni-\xi_n^2+4\xi_n^3i)}{(\xi_n-i)^5(\xi_n+i)^2}\sum_{j,l=1}^{n-1}X_jY_l\xi_j\xi_ld\xi_n\sigma(\xi')dx'\nonumber\\
&=-ih'(0)\Omega_3\int_{\Gamma^{+}}\frac{17\xi_ni-\xi_n^2+4\xi_n^3i}{(\xi_n-i)^5(\xi_n+i)^2}\sum_{j,l=1}^{n-1}X_jY_l\xi_j\xi_ld\xi_{n}dx'\nonumber\\
&=-ih'(0)\Omega_3\frac{2\pi i}{4!}\frac{4\pi}{3}\sum_{j}^{n-1}X_jY_j\left[\frac{17\xi_ni-\xi_n^2+4\xi_n^3i}{(\xi+i)^2}\right]^{(4)}\bigg|_{\xi_n=i}dx'\nonumber\\
&=-\frac{h'(0)}{2}\sum_{j=1}^{n-1}X_jY_j\pi^2\Omega_3dx'.
\end{align}
Let $X=X^T+X_n\partial_n,~Y=Y^T+Y_n\partial_n,$ then we have $\sum_{j=1}^{n-1}X_jY_j=g(X^T,Y^T).$ Now $\Phi$ is the sum of the cases (a), (b) and (c). Therefore, we get
\begin{align}\label{795}
\Phi=\sum_{i=1}^5\Phi_i=\left[\frac{13+32i}{32}X_nY_n+\frac{(13-10i)\pi}{24}g(X^T,Y^T)\right]h'(0)\pi\Omega_3dx'.
\end{align}
Then, by (\ref{b15})-(\ref{795}), we obtain following theorem
\begin{thm}\label{thmb1}
Let $M$ be a $4$-dimensional oriented
compact spin manifold with boundary $\partial M$ and the metric
$g^{M}$ be defined as (\ref{b1}), then we get the following equality:
\begin{align}
\label{b2773}
&\widetilde{{\rm Wres}}[\pi^+(\nabla_X^{S(TM)}\nabla_Y^{S(TM)}D^{-2})\circ\pi^+(D^{-2})]\nonumber\\
&=\frac{4\pi^2}{3}\int_{M}[Ric(X,Y)-\frac{1}{2}sg(X,Y)]d{\rm Vol_{M}}+\pi^2\int_{M}sg(X,Y)d{\rm Vol_{M}}\nonumber\\
&+\int_{\partial M}\left[\frac{13+32i}{32}X_nY_n+\frac{(13-10i)\pi}{24}g(X^T,Y^T)\right]h'(0)\pi\Omega_3d{\rm Vol_{M}}.
\end{align}
\end{thm}


\section*{Acknowledgements}
This work was supported by NSFC. 11771070 .
 The authors thank the referee for his (or her) careful reading and helpful comments.

\clearpage
\section*{Statement of ``The spectral Einstein functional and the noncommutative residue for manifolds with boundary"}
   a. Competing Interests: The authors have no relevant financial or non-financial interests to disclose.\\
\indent b.  Author Contribution Statement: All authors contributed to the study conception and design. Material preparation, data collection and analysis were performed by Tong Wu and Yong Wang. The first draft of the manuscript was written by Tong Wu and all authors commented on previous versions of the manuscript. All authors read and approved the final manuscript.\\
\indent c. Funding Information: This research was funded by National Natural Science Foundation of China: No.11771070.\\
\indent d. Availability of data and materials: The datasets supporting the conclusions of this article are included within the article and its additional files.\\
\end{document}